\documentclass[12pt,reqno,a4wide]{amsart}

\usepackage{tikz} 
\tikzset{
	arrow/.pic={\path[tips,every arrow/.try,->,>=#1] (0,0) -- +(.1pt,0);},
	pics/arrow/.default={triangle 45}
}
\usetikzlibrary{decorations.pathreplacing}
\usetikzlibrary{arrows,decorations.pathmorphing,snakes}
\tikzset{snake it/.style={decorate, decoration={snake, amplitude = 2, segment length = 5,}}}

\usetikzlibrary{decorations.pathreplacing}
\usetikzlibrary{arrows,decorations.pathmorphing,snakes}
\usetikzlibrary{fit}
\usepackage{calrsfs}
\usepackage[charter]{mathdesign}

\allowdisplaybreaks

\begin{document}
	
\title{Sums of powers over equally spaced {F}ibonacci numbers}

\author[Helmut Prodinger]{Helmut Prodinger}
\address{Department of Mathematics, University of Stellenbosch 7602,
	Stellenbosch, South Africa}
\email{hproding@sun.ac.za}


\begin{abstract} 
Recent results about sums of cubes of Fibonacci numbers ~\cite{Frontczak} are extended to arbitrary powers. 
\end{abstract}

\maketitle

\section{Introduction}

Frontczak~\cite{Frontczak} evaluates
\begin{equation*}
\sum_{k=0}^nF^3_{mk},\quad
\sum_{k=0}^n(-1)^kF^3_{mk},\quad
\sum_{k=0}^nL^3_{mk},\quad
\sum_{k=0}^n(-1)^kL^3_{mk},
\end{equation*}
with Fibonacci and Lucas numbers and $m$ being an odd integer.

We show here how to deal with general integer exponents (not just 3), and drop the restriction that $m$ must be odd.

Note that several papers about the evaluation of
\begin{equation*}
\sum_{k=0}^nF^{2l}_{2k}
\end{equation*}
and similar sums have been written a few years ago, here are just two citations \cite{ChuLi, Prodinger-Melham}. The instance where $2$ is replaced by $m$ is somewhat more delicate but quite instructive.

For later use, we mention the Binet formul\ae: Set
\begin{equation*}
\alpha=\frac{1+\sqrt5}{2}\quad\text{and}\quad \beta=\frac{1-\sqrt5}{2},
\end{equation*}
then $\sqrt5F_n=\alpha^n-\beta^n$ and $L_n=\alpha^n+\beta^n$.

\section{The summation of shifted Fibonacci numbers}

To deal with a sum over $F_{mn}$ (fixed $m$), we first consider a generating function:
\begin{equation*}
	\sum_{n\ge0}F_{nm}z^n=\frac{zF_m}{1-(F_{m+1}+F_{m-1})z+(-1)^mz^2}
	=\frac{zF_m}{1-L_mz+(-1)^mz^2}.
\end{equation*}
For a proof, rewrite it as
\begin{equation*}
	[1-(F_{m+1}+F_{m-1})z+(-1)^mz^2]\sum_{n\ge0}F_{nm}z^n=zF_m
\end{equation*}
and compare coefficients:
\begin{equation*}
	F_{nm}-(F_{m+1}+F_{m-1})F_{(n-1)m}+(-1)^mF_{(n-2)m}=[\!\![n=1]\!\!]F_m.
\end{equation*}
This can be proved by the Binet formula or otherwise and is  classical.

Consequently,
\begin{align*}
	&\frac1{1-z}\sum_{n\ge0}F_{nm}z^n\\
	&=\frac{F_m}{1-F_{m-1}-F_{m+1}+(-1)^m}\bigg[\frac1{1-z}-\frac{1-z(-1)^m}{1-(F_{m+1}+F_{m-1})z+(-1)^mz^2}\bigg].
\end{align*}	

Now we read off the coefficient of $z^n$:
\begin{align*}
\sum_{k=0}^nF_{km}	&=\frac{F_m}{1-F_{m-1}-F_{m+1}+(-1)^m}\\&-[z^n]\frac{F_m}{1-F_{m-1}-F_{m+1}+(-1)^m}\bigg[\frac{1-z(-1)^m}{1-(F_{m+1}+F_{m-1})z+(-1)^mz^2}\bigg]\\
	&=\frac{F_m}{1-F_{m-1}-F_{m+1}+(-1)^m}\\
	&-\frac{1}{1-F_{m-1}-F_{m+1}+(-1)^m}[z^n]\frac{F_m}{1-(F_{m+1}+F_{m-1})z+(-1)^mz^2}\\&+
	\frac{1}{1-F_{m-1}-F_{m+1}+(-1)^m}[z^n]\frac{-z(-1)^mF_m}{1-(F_{m+1}+F_{m-1})z+(-1)^mz^2} \\
	&=\frac{F_m-F_{(n+1)m}+(-1)^mF_{nm}}{1-F_{m-1}-F_{m+1}+(-1)^m}.
\end{align*}

So we notice that we have a closed formula for each fixed integer $m$.	

Let us also do the analogous computation for alternating sums:
\begin{equation*}
	\sum_{n\ge0}F_{nm}(-1)^nz^n=\frac{-zF_m}{1+(F_{m+1}+F_{m-1})z+(-1)^mz^2}
\end{equation*}
and
\begin{align*}
	\frac1{1-z}\sum_{n\ge0}F_{nm}(-1)^nz^n&=-{\frac 
		{F_m}{ 1+F_{m-1} +F_{m+1} + (-1) ^{m}   }}\frac1{1-z}\\*
	&+{\frac 
		{F_m}{  1+F_{m-1} +F_{m+1} + (-1) ^{m}   }}{\frac { 1-(-1)^mz}{    1+z(F_{m-1} +F_{m+1}) + (-1) ^{m}{z}^{2}   }}.
\end{align*}

Reading off the coefficient of $z^n$ on both sides leads to
\begin{align*}
	\sum_{k=0}^n(-1)^kF_{mk}=\frac 
	{-F_m+(-1)^nF_{(n+1)m}+(-1)^{n+m}F_{nm}}{ 1+F_{m-1} +F_{m+1} + (-1) ^{m}   }.
\end{align*}

\section{Summing shifted Lucas numbers}

First, we need the generating function
\begin{equation*}
	\frac{2-zL_m}{1-zL_m+z^2(-1)^m}=\sum_{k\ge0}L_{mk}z^k,
\end{equation*}
which holds for $m\ge0$.

It can again be checked by writing it as
\begin{equation*}
	2-zL_m=(1-zL_m+z^2(-1)^m)\sum_{k\ge0}L_{mk}z^k,
\end{equation*}
comparing coefficients and prove that the resulting coefficients are zero for $k\ge2$, either by the Binet formula or
by using classical identities for Lucas numbers.

Furthermore, for $m\ge1$,
\begin{equation*}
	\frac1{1-z}\frac{2-zL_m}{1-zL_m+z^2(-1)^m}=
	\frac{1-L_m}{1-L_m+(-1)^m}\frac1{1-z} +\frac{1}{1-L_m+(-1)^m}\frac{(-1)^m(1+z-zL_m)}{1-zL_m+z^2(-1)^m}.
\end{equation*}
Comparing coefficients of $z^n$, this leads to
\begin{align*}
	\sum_{k=0}^{n}L_{mk}&=\frac{1-L_m}{1-L_m+(-1)^m}\\&+
	\frac{(-1)^m}{1-L_m+(-1)^m}[z^n]\frac{1}{1-zL_m+z^2(-1)^m}+
	\frac{(-1)^m(1-L_m)}{1-L_m+(-1)^m}[z^n]\frac{z}{1-zL_m+z^2(-1)^m}\\
	&=\frac{1-L_m}{1-L_m+(-1)^m}+
	\frac{(-1)^m}{1-L_m+(-1)^m}\frac{F_{(n+1)m}}{F_m}+
	\frac{(-1)^m(1-L_m)}{1-L_m+(-1)^m}\frac{F_{nm}}{F_m}.
\end{align*}

We can deal with an alternating version by small modifications:
\begin{equation*}
	\frac{2+zL_m}{1+zL_m+z^2(-1)^m}=\sum_{k\ge0}(-1)^kL_{mk}z^k,
\end{equation*}
and
\begin{align*}
	\frac1{1-z}	\frac{2+zL_m}{1+zL_m+z^2(-1)^m}&=\frac{2+L_m}{1+L_m+(-1)^m}\frac1{1-z}\\
	&+\frac{2(-1)^m+L_m+2z(-1)^m+zL_m(-1)^m}{1+L_m+(-1)^m}\frac{1}{1+zL_m+z^2(-1)^m}.
\end{align*}
Reading off coefficients,
\begin{align*}
	\sum_{k=0}^n(-1)^kL_{mk}&=\frac{2+L_m}{1+L_m+(-1)^m}\\
	&+\frac{2(-1)^m+L_m}{1+L_m+(-1)^m}(-1)^n\frac{F_{(n+1)m}}{F_m}+\frac{2(-1)^m+L_m(-1)^m}{1+L_m+(-1)^m}(-1)^{n-1}\frac{F_{nm}}{F_m}.
\end{align*}

\section{Expanding powers of Fibonacci and Lucas numbers}

Our goal is  here to expand $F_n^j$ in terms of $F_{mn}$, and likewise for Lucas numbers. To clarify, we start with a list of such expansions:

\begin{align*}
	F_n^2&=\frac25F_{2(n+1)}-\frac35F_{2n}-\frac25(-1)^n\\
	F_n^3&=\frac15F_{3n}-\frac35(-1)^nF_n\\
	F_n^4&=\frac{2}{75}F_{4(n+1)}-\frac{7}{75}F_{4n}-\frac{8}{25}(-1)^nF_{2(n+1)}+\frac{12}{25}(-1)^nF_{2n}+\frac{6}{25}\\
	F_n^5&=\frac1{25}F_{5n}-\frac15(-1)^nF_{3n}+\frac25F_n\\
	F_n^6&=\frac1{500}F_{6(n+1)}-\frac9{500}F_{6n}-\frac4{125}(-1)^nF_{4(n+1)}+\frac{14}{125}(-1)^nF_{4n}\\&\qquad+\frac6{25}F_{2(n+1)}-\frac9{25}F_{2n}-\frac4{25}(-1)^n\\
	F_n^7&=\frac1{125}F_{7n}-\frac7{125}(-1)^nF_{5n}+\frac{21}{125}F_{3n}-\frac{7}{25}(-1)^nF_n\\
\end{align*}

The formula for $F_n^j$, odd $j$, is easier to guess:
\begin{equation*}
F_n^j=\frac1{5^{(j-1)/2}}\sum_{0\le s<\frac j2}F_{(j-2s)n}\binom{j}{s}(-1)^{sn}.
\end{equation*}
The instance $j$ being odd is harder, but here is the result:
\begin{align*}
F_n^j&=\frac1{5^{j/2}}\sum_{1\le s\le j/2}\frac{2}{F_{2s}}(-1)^{(n+1)(\frac j2+s)}\binom{j}{j/2+s}F_{2s(n+1)}\\
&-\frac1{5^{j/2}}\sum_{1\le s\le j/2}\frac{L_{2s}}{F_{2s}}(-1)^{(n+1)(\frac j2+s)}\binom{j}{j/2+s}F_{2sn}\\
&+\frac1{2\cdot5^{j/2}}\binom{j}{j/2}\big(1-(-1)^n+(-1)^{j/2}+(-1)^{n+j/2}\big).
\end{align*}

	The results for Lucas numbers are somewhat simpler:
	\begin{align*}
L_n^2&=L_{2n}+2(-1)^n\\
L_n^3&=L_{3n}+3(-1)^nL_{n}\\
L_n^4&=L_{4n}+4(-1)^nL_{2n}+6\\
L_n^5&=L_{5n}+5(-1)^nL_{3n}+10L_{n}\\
L_n^6&=L_{6n}+6(-1)^nL_{4n}+15L_{2n}+20(-1)^n\\
L_n^7&=L_{7n}+7(-1)^nL_{5n}+21L_{3n}+35(-1)^nL_{n}\\
L_n^8&=L_{8n}+8(-1)^nL_{6n}+28L_{4n}+56(-1)^nL_{2n}+70\\
	\end{align*}
It is not too hard to guess the general formula from that:
\begin{equation*}
L_n^j=\sum_{0\le s <\frac j2}L_{n(j-2s)}\binom{j}{s}(-1)^{sn}+(-1)^n[\!\![j\text{ even}]\!\!]\binom{j}{j/2}.
\end{equation*}	
	
	Once these formul\ae\ have been successfully guessed (the hard part), they can be proved using the Binet
	formul\ae\ and routine manipulations with binomial identities. We leave this for the interested reader.

	Summations like
	\begin{equation*}
		\sum_{0\le k\le n}1=n+1\qquad\text{and}\qquad
		\sum_{0\le k\le n}(-1)^k=\frac12\big(1+(-1)^n\big)
	\end{equation*}
	are also needed but of a trivial nature.
	
	\section{Frontczak's results revisited}
	
	Let us do an example computation:
	\begin{align*}
\sum_{k=0}^nF_{mk}^3&=
\sum_{k=0}^n\Big[\frac15F_{3mk}-\frac35(-1)^kF_{mk}\Big]\\
&=\frac15	\frac{F_{2m}-F_{(n+1)3m}+(-1)^mF_{3nm}}{1-F_{3m-1}-F_{3m+1}+(-1)^m}
-\frac35	\frac{-F_m+(-1)^nF_{(n+1)m}+(-1)^{n+m}F_{nm}}{1-F_{m-1}-F_{m+1}+(-1)^m}.
	\end{align*}
The other sums from \cite{Frontczak} can be obtained in a similar way.
	
\section{Why can we expand powers of Fibonacci and Lucas numbers?}	

The key to the success is the formula
\begin{equation*}
x^n+y^n=\sum_{0\le k \le n/2}(-1)^k\frac{n}{n-k}\binom{n-k}{k}(x+y)^{n-2k}(xy)^k,
\end{equation*}
which is a consequence of classical formul\ae\ due to Girard and Waring; see e. g. \cite{Gould}.

Set
\begin{equation*}
\alpha=\frac{1+\sqrt5}{2}\quad\text{and}\quad \beta=\frac{1-\sqrt5}{2},
\end{equation*}
then, with $x=\alpha^m$ and $y=\beta^m$, the formula becomes
\begin{equation*}
	L_{mn}=\sum_{0\le k \le n/2}\frac{n}{n-k}\binom{n-k}{k}L_{(n-2k)m};
\end{equation*}
if one sets $x=\alpha^m$ and $y=(-\beta)^m$, then
\begin{equation*}
x+y=\begin{cases}L_m &\text{if $m$ is even},\\
	\sqrt{5}F_m&\text	{if $m$ is odd}
	\end{cases}
\end{equation*}
and
\begin{equation*}
	x^n+y^n=\begin{cases}L_{mn} &\text{if $m$ or $n$ is even},\\
		\sqrt{5}F_{mn}&\text	{if $m$ and $n$ are odd}.
	\end{cases}
\end{equation*}
Furthermore, $xy=1$. Since $L_n=F_{n+1}+F_{n-1}$, everything could be expressed in Fibonacci number (alternatively, everything could be expressed in terms of Lucas numbers).

In the paper \cite{Prodinger-Melham}, these formul\ae\ were derived from scratch.

So, $F_{mn}$ (resp.\ $L_{mn}$) are expressed in terms (linear combinations) of powers of $F_m$ resp.\ $L_m$. The formul\ae\ 
of the previous sections are just inverted versions of this, namely, powers of Fibonacci numbers are expressed as linear combinations
of shifted Fibonacci numbers.

\bibliographystyle{plain}


\end{document}